\pdfoutput=1
\documentclass[onecolumn,12pt]{article}

\usepackage{arxiv}

\newcommand*{\rootPath}{.}
\usepackage{parskip}

\usepackage[small,labelfont={bf,sf}]{caption}

\usepackage{abstract}

\usepackage{hyperref} %deactivate for arXiv upload
\hypersetup{colorlinks,
            citecolor=mdarkgreen,
            linkcolor=mred}
\AtBeginDocument{}%
\AtBeginDocument{}%

\usepackage{amsmath}
\usepackage{amssymb}
\usepackage{amsthm}
\usepackage{mathtools}
\usepackage[ruled,vlined, linesnumbered]{algorithm2e}
\usepackage{siunitx}
\usepackage{xargs}
\usepackage[table]{xcolor}  % Coloured text etc.
\usepackage{enumerate}
\usepackage{balance}

\usepackage{csquotes} % for proper quotations
\usepackage{csvsimple} % to generate the benchmark result table
\usepackage{booktabs} % to generate the benchmark result table
\usepackage[para,online,flushleft]{threeparttable}
\usepackage{multirow}

\usepackage{tikz}
\usepackage{pgfplots}
\usepackage{gincltex}
\usepackage{pgfmath,pgffor}
\usepackage{xargs}                      % Use more than one optional parameter in a new commands
\pgfplotsset{compat=newest} %globally avoid intereference of ticks and yaxis labels in plots with tikz package
\pgfplotsset{
      table/search path={\rootPath/figs/data},
    }
\usetikzlibrary{shapes.geometric,calc,arrows,intersections}

\newcommand{\norm}[1]{\left\lVert#1\right\rVert}

\newcommand{\mc}[1]{\mathcal{#1}}
\newcommand{\mb}[1]{\mathbb{#1}}

\DeclareMathOperator*{\argmin}{argmin}

\theoremstyle{definition}

\theoremstyle{remark}

\makeatletter
\newcommand{\removelatexerror}{\let\@latex@error\@gobble}
\makeatother

\definecolor{mblue}{HTML}{1F77B4}
\definecolor{morange}{HTML}{FF7F0E}
\definecolor{mred}{HTML}{D62728}
\definecolor{mpurple}{HTML}{9467BD}
\definecolor{mbrown}{HTML}{8C564B}
\definecolor{mgreen}{HTML}{00802b}

\definecolor{mmred}{HTML}{e0024c}
\definecolor{mdarkgreen}{RGB}{11, 91, 35}

\definecolor{mosekcolor}{rgb}{0.863,0.129,0.302}% or {220,33,77} on a 255 scale
\definecolor{scscolor}{rgb}{0,0.549,0}   % or {0,140,0} on a 255 scale
\definecolor{cosmocolor}{rgb}{0,0.392,0.871} % or {0,148,222} on a 255 scale
\definecolor{cosmocdcolor}{rgb}{1,0.4,0} % or {102,0,102} on a 255 scale
\definecolor{osqpcolor}{rgb}{0.4,0.0,0.4} % or {55,200,171} on a 255 scale
\definecolor{gurobicolor}{RGB}{11, 91, 35}

\newcommand{\pkg}[1]{\texttt{#1}}

\DeclareMathOperator{\prox}{prox}

\renewcommand{\Re}{\mathbb{R}}
\newcommand{\eqdef}{\coloneqq}

\newcommand{\vacc}{v_k^{\text{acc}}}
\let\OLDthebibliography\thebibliography
\renewcommand\thebibliography[1]{
  \OLDthebibliography{#1}
  \setlength{\parskip}{0pt}
  \setlength{\itemsep}{0pt plus 0.3ex}
}

\usepackage{graphicx}
\title{\LARGE \bf Safeguarded Anderson acceleration for parametric nonexpansive operators}

\author{Michael Garstka, Mark Cannon \& Paul Goulart$^{1}$% <-this % stops a space
\thanks{$^{1}$The authors are with the Department of Engineering Science, University of Oxford, Oxford, OX1 3PJ, UK. Email: {\tt\{michael.garstka, mark.cannon, paul.goulart\}@eng.ox.ac.uk}}
}

\begin{document}
\maketitle
\thispagestyle{plain}
\pagestyle{plain}

\begin{abstract}
This paper describes the design of a safeguarding scheme for Anderson acceleration to improve its practical performance and stability when used for first-order optimisation methods.
We show how the combination of a non-expansiveness condition, conditioning constraints, and memory restarts integrate well with solver algorithms that can be represented as fixed point operators with dynamically varying parameters.
The performance of the scheme is demonstrated on seven different QP and SDP problem types, including more than 500 problems.
%Our implementation of safeguarded Anderson acceleration is the default in the recent version of our
The safeguarded Anderson acceleration scheme proposed in this paper is implemented in the
open-source ADMM-based conic solver \pkg{COSMO}.
\end{abstract}
\section{INTRODUCTION}\label{sec:introduction}
\noindent Solutions of large convex optimisation problems of the form
\begin{equation}
  \begin{array}{ll}
    \mbox{minimize}   & \textstyle{\frac{1}{2}}x^\top Px + q^\top x\\
    \mbox{subject to} & Ax + s  = b, \, (x, s) \in \Re^n \times \mc{K}
  \end{array}
  \label{eq:primal}
\end{equation}
with $P\in \mathbb{S}^n_+, \, q \in \Re^n, \, A \in \Re^{m \times n}, \, b \in \Re^m$,
and $\mc{K}$ a convex cone, are vital to numerous application areas, including robust and optimal control, structural design, operations research, and signal processing~\cite{Boyd_2004,Wolkowicz_2012,BenTal_2001,Mattingley_2010}.
Moreover, recent interest in machine learning has seen convex optimisation being used in many novel applications, including neural network verification against adversarial attacks~\cite{Raghunathan_2018}, sparse principal component analysis~\cite{d_aspremont_2004}, graph clustering~\cite{Bie_2006}, and kernel matrix learning~\cite{Lanckriet_2004}.
Many of these applications require the solution of very large-scale problem instances, which presents a challenge for established interior-point solution methods (IPMs).
This is because IPMs solve a Newton system at each iteration that grows with the problem dimension~$n$ and has a per iteration computational cost of $\mc{O}(n^3)$.
This drawback led to a renewed interest in first-order methods (FOMs) that trade-off moderate accuracy solutions for a lower per-iteration computational cost, allowing them to solve larger problems.
Popular FOMs  that have been developed into solver packages include the \emph{Alternating Direction Method of Multipliers (ADMM)}~\cite{Stellato_2020, Zheng_2020, Garstka_2020}, \emph{Douglas-Rachford splitting}~\cite{ODonoghue_2016}, and the \emph{Augmented Lagrangian method}~\cite{Zhao_2010}.
Moreover, the \emph{proximal gradient method}~\cite{Beck_2009} is a popular tool to design custom solution algorithms.

The minimisation problem solved by FOMs can be recast as the problem of finding a fixed point of a nonexpansive operator $F \colon \mc{D} \subseteq \Re^n \rightarrow \Re^n$
\begin{equation}\label{eq:fixed_point_iter}
  v = F(v),
\end{equation}
where $v \in \Re^n$.
%For the following discussion of different acceleration methods, we assume that they are applied to a FOM that can be represented in the form~\eqref{eq:fixed_point_iter} and that conditions such as non-expansiveness hold as required.
Several recent publications consider acceleration methods for FOMs (see the survey in~\cite{dAspremont_2021}).
Well-known acceleration methods based on the gradient descent method include \emph{Nesterov's accelerated gradient method}, \emph{accelerated proximal gradient method} and the \emph{Heavy ball method}.
However, all of these methods require differentiability of at least part of the objective of the underlying optimisation problem.

A line search method for general nonexpansive operators is developed in~\cite{Giselsson_2016}.
This searches in the direction of the fixed-point residual
%along the last residual direction
and takes larger steps than classical line search methods.
Although it is possible to skip many iterations with this approach, it requires the operator $F$ to be evaluated at each candidate point, incurring the same cost as a full iteration.
Therefore, this method is only practical if the operator can be expressed in the form $F(v) = F_2\bigl( F_1(v) \bigr)$ where $F_2$ is cheap to evaluate and $F_1$ is affine, which allows many points in the search direction to be checked cheaply.

Many acceleration approaches view~\eqref{eq:fixed_point_iter} as the problem of finding the zeros $0 \in r(v) $ of the residual operator $r(v)= F(v) - v$, and apply Newton's method due to its fast practical asymptotic convergence.
For example, in~\cite{Ali_2017} a semismooth Newton method is applied to the fixed-point operator that is used in the \pkg{SCS}~\cite{ODonoghue_2016} solver.
The method assumes a semismooth operator and relies on expensive line search steps.
Similarly, the \emph{SuperMann} acceleration scheme~\cite{Themelis2019} globalizes fixed-point iterations of non-expansive operators and enjoys superlinear convergence under certain conditions.
This approach uses a limited memory Broyden method to approximate the Jacobian with Powell regularisation to keep the Jacobian approximation non-singular.
The original fixed point iteration is employed if the accelerated candidate point does not provide an improved solution estimate.

%A number of recent publications have suggested a similar approach, but instead of using Newton’s or Broyden’s method the update equation for the approximate Jacobian is based on \emph{Anderson acceleration} (AA)~\cite{Anderson_1965}, which does not require a line search.
A similar approach uses \emph{Anderson acceleration} (AA)~\cite{Anderson_1965}, which does not require a line search, to update the  Jacobian  approximation instead of using Newton’s or Broyden’s method.
The method has been applied to electronic structure computation, known as \emph{Pulay mixing}~\cite{Pulay_1980}, \emph{direct inversion in the iterative subspace (DIIS)}~\cite{Pulay_1982}, and \emph{Anderson mixing}, but it has only recently gained attention in the optimisation community
(e.g.~\cite{Walker_2011,Fang_2009}).

The AA algorithm uses a combination of past iterates to find the accelerated point, where the weights are determined by minimizing a weighted sum of past residuals.
Unfortunately, convergence guarantees for AA are available only if additional assumptions such as contractivity~\cite{Toth_2015}, linearity~\cite{Potra_2013}, or differentiability~\cite{Gay_1978} are made on the operator or on the memory-length of the scheme.
However, FOMs for convex conic optimisation problems have none of these properties, and for large problems only a limited memory variant of AA is economical. Moreover, Mai and Johansson~\cite{Mai_2020} prove that the limited-memory AA cannot guarantee global convergence.
Thus, additional safeguarding measures are required.
An implementation to safeguard AA in combination with the Douglas-Rachford method employed by \pkg{SCS} is discussed in~\cite{Zhang_2020}.
This combines (type-I) AA steps with the execution of the original algorithm whenever the residual decreases sufficiently.
In order to ensure a non-singular Jacobian, a Gram-Schmidt orthogonalization strategy is used, and a \emph{rolling-memory} approach is also employed.

%\medskip
\textit{Contributions of this paper:}
\begin{enumerate}[1.]
  \item We design a safeguarding mechanism for AA based on a relaxed non-expansiveness condition on the norm of the residual operator.
  \item We investigate a restarted limited-memory scheme that allows the use of a solver algorithm whose representation is a parametric operator and which relies periodically on non-accelerated steps for infeasibility detection.
  \item We provide evidence from more than 500 QP and SDP test problems demonstrating the performance of our implementation against the non-accelerated operator.
The mean number of iterations for each problem set is reduced by a factor between \num{1.7} and \num{8.5} and the mean solve time by a factor of up to $6$ for higher accuracy solutions.
We also show that the additional time required to calculate the Anderson directions can be kept between 3\% to 15\% for large QPs and SDPs.
  Thus, the results make a particularly strong case for using AA to solve SDPs and large QPs.
\end{enumerate}

%\medskip
\textit{Outline:}
In Section~\ref{sec:background} we review the relationship between solving a convex conic problem with a FOM and applying a fixed point iteration to a particular operator.
Then we introduce the classical AA method and explain how it can be viewed as a multisecant Broyden's method.
In Section~\ref{sec:safeguarded} we detail the design decisions of our safeguarding scheme for AA.
Section~\ref{sec:benchmarks} shows benchmark results of different variants of the \pkg{COSMO} solver: the classic algorithm, the accelerated algorithm, and the safeguarded \emph{and} accelerated algorithm.
Section~\ref{sec:conclusion} concludes the paper.

%\medskip
\textit{Notation:}
Denote the space of real numbers $\mb{R}$, the $n$-dimensional real space $\mb{R}^n$, and denote a convex proper cone as $\mc{K}$.
The identity operator is denoted as $\text{Id}$.
We say that an operator $F \colon \mc{D} \subseteq \Re^n \rightarrow \Re^n$ is:
\begin{enumerate}[a)]
  \item \emph{nonexpansive} if $\norm{Fx  - Fy} \leq \norm{x - y}$ for all $x,y \in \Re^n$;
  \item \emph{$\alpha$-averaged} if there exists a nonexpansive operator $G\colon \Re^n \rightarrow \Re^n$ such that $F = (1-\alpha)\text{Id} + \alpha G$ for $\alpha \in (0,1)$;
  \item \emph{firmly nonexpansive} if it is $\frac{1}{2}$-averaged.
\end{enumerate}

The scaled \emph{proximal operator} of a convex, closed and proper function $f\colon \mb{R}^n \rightarrow \mb{R} \cup \{\infty\}$ with $\gamma > 0$ is given by
\begin{equation}\label{eq:prox}
  \prox_{\gamma f}(v) \eqdef \argmin_y \{ f(y) + \textstyle{\frac{1}{2\gamma}}\norm{y-v}^2_2 \}.
\end{equation}
Denote the \emph{reflected proximal operator} as $R_{\gamma f}(v) \coloneqq  (2 \prox_{\gamma f} - \text{Id})(v)$.
We note that the proximal operator \eqref{eq:prox} is firmly nonexpansive and the reflected proximal operator is nonexpansive~\cite{Bauschke_2011}.

\section{BACKGROUND}
\label{sec:background}
Before introducing the AA method, we review the relationship between solving an optimisation problem via a FOM and finding the fixed points of a firmly nonexpansive operator.

\subsection{FOMs as fixed-point iterations}

In this paper we propose a safeguarded AA method to improve the slow convergence of FOMs when solving large convex optimisation problems of the form~\eqref{eq:primal}.
A popular FOM to solve~\eqref{eq:primal} to moderate accuracy is \emph{Douglas-Rachford splitting (DRS)}. This solves problems of the form
\begin{equation}
  \begin{array}{ll}
    \mbox{minimize}   & f(z) + g(z),
  \end{array}
  \label{eq:conic_prob}
\end{equation}
where both $f \colon \Re^n \rightarrow \Re \cup \{\infty\}$ and $g\colon \Re^n \rightarrow \Re \cup \{\infty\}$ are convex, closed, and proper.
It is well known that the Fenchel dual of~\eqref{eq:primal} is in the form of~\eqref{eq:conic_prob}. Applying DRS to the problem gives the following algorithm
\begin{subequations}
  % \noeqref{eq:dr_z,eq:dr_x, eq:dr_v}
\begin{align}
  z^k &\coloneqq \mathbf{prox}_{\gamma f}(v^k) \label{eq:dr_z}\\
  x^k &\coloneqq \mathbf{prox}_{\gamma g}(2 z^k - v^k) \label{eq:dr_x}\\
  v^{k+1} &\coloneqq v^k + (x^k - z^k)\label{eq:dr_v}
\end{align}
\end{subequations}
with $\gamma > 0$, or using the more compact operator form
\begin{subequations}
\label{eq:DRexp}
\begin{align}
v^{k+1} &\coloneqq F(v^k) = \left(\frac{1}{2}\text{Id} + \frac{1}{2} R_{\gamma f} R_{\gamma g}\right)(v^k), \\
z^{k+1} & \coloneqq \mathbf{prox}_{\gamma f}(v^{k+1}).
\end{align}
\end{subequations}
Notice that~(\ref{eq:DRexp}a,b) has the fixed-point operator form~\eqref{eq:fixed_point_iter}.
Moreover $F$ is a $\frac{1}{2}$-averaged iteration of the nonexpansive operator $ R_{ \gamma f} R_{ \gamma g}$ and therefore firmly-nonexpansive.
Hence applying the \emph{Picard iteration}~\eqref{eq:fixed_point_iter} will converge linearly to a fixed point (if one exists) that is coincident with the optimal solution of the underlying optimisation problem~\cite{Ryu_2016}.

% We further assume that to influence the convergence behaviour of the operator it be changed using a parameter $\rho$, denoted as $F_\rho$. While changes to $\rho$ might happen multiple times as we iterate we also assume that after a while we keep $\rho$ constant.
% For every $\rho$ (?), $F_\rho$ is firmly-nonexpansive and the Picard-iteration~\eqref{eq:fixed_point_iter} will converge to a fixed point if one exists~\cite{Ryu_2016}.
%For the remainder of this paper we view
We consider the fixed-point operator $F$ in (\ref{eq:DRexp}a)  as an abstraction of the actual FOM solver used to solve~\eqref{eq:primal}.

\subsection{Anderson acceleration}
AA calculates an accelerated candidate point $\vacc$ as a weighted combination of $m_k+1$ previous iterates,
\begin{equation}\label{eq:anderson_update}
  \vacc = \sum_{i=0}^{m_k} \alpha_k^i F(v_{k-m_k + i}).
\end{equation}
Anderson's main idea was to choose the weights $\alpha_k = \begin{bmatrix}\alpha_k^0, \dots, \alpha_k^{m_k} \end{bmatrix}$ by minimizing the norm of past residual vectors:
\begin{equation}\label{eq:anderson_alphas}
  \begin{array}{ll}
    \mbox{minimize}   & \norm{R_k \alpha_k }_2^2\\
    \mbox{subject to} & \mathbf{1}_{m_k+1}^\top \alpha_k = 1,
  \end{array}
\end{equation}
where $R_k = \begin{bmatrix}
r_{k-m_k} & \cdots & r_k
\end{bmatrix}$ is the matrix of past residual vectors, $r_k \coloneqq r(v_k)$.
Eyert~\cite{Eyert_1996} established a connection between AA and multisecant Broyden's method which allows AA to be considered a Quasi-Newton method.
To clarify the relationship one can make a change of variables~\cite{Fang_2009} by defining $\eta = (\eta^0, \ldots, \eta^{m_k-1}) \in \Re^{m_k}$, which implicitly encodes the constraint in~\eqref{eq:anderson_alphas}
\begin{equation}
  \alpha^0 = \eta^0,\, \alpha^i = \eta^i - \eta^{i-1}, \ldots, \, \alpha^{m_k} = 1 - \eta^{m_k - 1}.
\end{equation}
Next, $\eta$ is substituted for $\alpha$ and we use the differences between iterates $\Delta v_k = v_{k+1}- v_{k}$, $\Delta r_k = r(v_{k+1}) - r(v_{k})$ to rewrite~\eqref{eq:anderson_update} and~\eqref{eq:anderson_alphas} as
\begin{align}
  \vacc &= v_k - r_k  - (\mc{V}_k - \mc{R}_k)\eta_k,
\end{align}
$\mc{V}_k = [ \Delta v_{k-m_k}, \ldots, \Delta v_{k-1} ]$, $\mc{R}_k = [ \Delta r_{k-m_k}, \ldots, \Delta r_{k-1} ]$.
The Anderson coefficients $\eta$ are calculated by solving the least-squares problem $\norm{r_k - \mc{R}_k \eta}_2^2$.
Assuming $\mc{R}_k$ is full-rank we can substitute the solution for $\eta$ and write AA as a multisecant Broyden type-II method:
\begin{align}
  \vacc &= v_k - H_k^{II} r_k \label{eq:vacc_update}
\end{align}
where $H_k^{II} =  I +  (\mc{V}_k - \mc{R}_k)(\mc{R}_k^\top \mc{R}_k)^{-1}\mc{R}_k^\top$ can be viewed as the rank-$m_k$ update formula to approximate the inverse of the Jacobian of the residual operator.
A corresponding type-I method is obtained by using the approximation update rule $H_k^I =  I +  (\mc{V}_k - \mc{R}_k)(\mc{V}_k^\top \mc{R}_k)^{-1}\mc{V}_k^\top$.

\section{SAFEGUARDED ANDERSON ACCELERATION}\label{sec:safeguarded}
We make no assumptions about the smoothness of the operator $F$, which means global convergence of the unaltered acceleration method cannot be guaranteed.
Section~\ref{sec:background} explains how Anderson acceleration can be seen as a Broyden's method which starts with $I$ as an estimate for the Jacobian inverse $H_k$ of the operator and then uses rank-$1$ updates to modify it based on input and output differences of the residual operator. A source of instability of this Jacobian update step is the tendency of FOMs to eventually reach regimes of slow convergence.
This causes the vectors used in the updates to become nearly collinear and leads to poor conditioning of the Jacobian approximation.
%The other weakness of AA, inherited from the class of Quasi-Newton methods, is the lack of global convergence guarantees.
Like other Quasi-Newton methods, AA also suffers from a lack of global convergence guarantees.
Therefore further safeguarding measures are needed, usually to check the quality of candidate iterates before they are accepted.

Solvers based on FOMs rarely stick to the classic algorithm but instead use heuristic rules that adapt some of the algorithm's parameters. Thus, the corresponding operator will change at different points in the solution process.
We denote the parameter-dependent operator $F_\rho$ with changing parameter vector $\rho$.
Our AA scheme needs to be able to accommodate changing operators.
The following subsections detail the design choices that are summarized in Algorithm~\ref{alg:safeguarded_acceleration}.

\begin{figure}[htb]
%\vspace*{0.15cm}
 \removelatexerror
  \begin{algorithm}[H]
\SetKwInOut{Input}{Input}
\Input{ $v_0, f_0$, fixed-point iteration $F_{\rho} \colon \Re^n \rightarrow \Re^n$ with $v_{k+1}=f_k = F_\rho(v_k)$, allocated memory for $\mc{V}_0$, $\mc{R}_0$, $Q$, and $R$,  column pointer $j = 1$, parameters: $ \eta_
{\text{max}}, \tau, \epsilon, m_{\text{max}}$}
 \texttt{acc\_success} = false\;
 \While{$\norm{v_{k+1} - v_k } > \epsilon$}{
 Update history: $\mc{V}_{j} \leftarrow [\mc{V}_{j-1}, \Delta v_{k-1}],$ $\mc{R}_{j} \leftarrow [\mc{R}_{j-1}, \Delta r_{k-1}]$,
 $j \leftarrow j + 1$\;
 \If{$j > 2$}{
 Update QR factors $Q_k$, $R_k$ of $\mc{R}_k$\;
 Compute $\eta_k$ from $R_k \eta_k = Q_k^\top r_k $\;
\eIf{ $\norm{\eta_k}_2 > \eta_{\text{max}}$}{
  \texttt{acc\_success} $\leftarrow$ false\;
}{
 Accelerate: $v_k^{acc} =  f_k - (\mc{V}_j -  \mc{R}_j)\eta_k$\;
 \texttt{acc\_success} $\leftarrow$ true\;
 }
 \If{\texttt{acc\_success}}{
    Fixed-point iteration: $f_{k}^{acc} = F_{\rho}(v_k^{acc})$\;
    Compute residual: $r_{k}^{acc} = v_k^{acc} - f_{k}^{acc}$\;
    \eIf{$\norm{r_{k}^{acc}}_2 \leq \tau \norm{r(v_{k-1})}_2$}{
    $v_{k+1} \leftarrow v_k^{acc}$, $f_{k+1} \leftarrow f_k^{acc}$, and $r_{k+1} \leftarrow r_k^{acc}$\;
    }{
   \texttt{acc\_success} $\leftarrow$ false\;
    }
}
} %j>2
\If{not \texttt{acc\_success} or $j \leq 2$}{
     \If{operator change scheduled}{
        $\rho^{*} \leftarrow u(F_{\rho}, v_k, f_k, r_k, \rho)$ \;
        Update operator: $F_\rho \leftarrow F_{\rho^*}$\;
        $j \leftarrow 1$\;
      }
  Safeguarding step: $v_{k+1} \leftarrow f_k$, $f_{k+1} = F_{\rho}(v_{k+1})$, $r_{k+1} = v_{k+1} - f_{k+1}$\;
}
    \If{infeasibility detection scheduled and $j = 2$}{
      perform infeasibility checks\;
    }
  \If{$j > m_{\text{max}}$}{
  $j \leftarrow 1$ \;
  }
 } %while
  \caption{Safeguarded AA with memory restarts and scheduling.}
   \label{alg:safeguarded_acceleration}
\end{algorithm}
\vspace{-4mm}
\end{figure}

\subsection{Anderson acceleration variant}
As discussed in Section~\ref{sec:background} and shown in~\cite{Fang_2009} two types of AA can be derived from the corresponding Broyden's methods.
These differ in the way that the Anderson coefficients $\eta$ are computed.
The type-I variant uses
%results in the equation
$\eta_k=(\mc{V}_k^\top \mc{R}_k)^{-1}\mc{V}_k^\top r_k$ whereas type-II uses $\eta_k = (\mc{R}_k^\top \mc{R}_k)^{-1}\mc{R}_k^\top r_k$.
For neither of the two variants we observed consistently faster convergence when benchmarked against standard problem sets~\ref{sec:benchmarks}.
However, the type-II variant has the advantage that the computation of $\eta_k$
%is in the form of a least-squares problem and can thus
can be solved using QR decomposition.
Moreover, the QR factorisation $\mc{R}_k = Q_kR_k$ can be efficiently updated from $Q_{k-1}, R_{k-1}$ as at each step only one column is added to $\mc{R}_k$~\cite{Walker_2011}.
Compared to the type-I computation of the coefficients, this results in a significant speed-up.

\subsection{Safeguarding mechanism}
One source of instability of AA is that not every acceleration candidate point is guaranteed to be closer to the fixed-point than the current iterate.
To mitigate this issue, we perform a safeguarding step at each iteration.
Assume that at iteration $k$, AA produced an accelerated candidate point $\vacc$ using~\eqref{eq:vacc_update}.
An intuitive way to assess the quality of $\vacc$ is by comparing the resulting residual operator norm with the last accepted step and imposing the condition
\begin{equation}\label{eq:safeguarding_check}
  \norm{r(v_k^\text{acc})}_2 = \norm{v_k^\text{acc} - F_\rho(v_k^\text{acc}) }_2 \leq \tau \norm{r(v_{k})}_2 ,
\end{equation}
where $\tau \in (0,1)$ is an expansiveness tolerance.
If~\eqref{eq:safeguarding_check} holds, then $v_k^{\text{acc}}$ is accepted as the next iterate, $v_{k+1} = v_k^{\text{acc}}$.
Otherwise, the method resorts to $v_{k+1} = F_\rho(v_k)$.
Checking condition~\eqref{eq:safeguarding_check} is expensive as it involves both the evaluation of $F_\rho(\vacc)$ and $F_\rho(v_k)$, which means two iterations of the underlying FOM.
This makes the condition inefficient in practice.
We also observed from standard benchmark testing that enforcing strict monotonicity of the residual norm tends to make the performance of the solver worse.
The authors in~\cite{Themelis2019} therefore only enforce the condition periodically or when past progress has stalled.
On the other hand, the term $\tau \norm{r(v_k)}_2$ is replaced in~\cite{Fu_2020} by a summable and exponentially decaying series based on $\norm{r(v_0)}_2$.

We use the relaxed safeguarding condition
\begin{equation}\label{eq:relaxed_safeguarding}
   \norm{v_k^\text{acc} - F_\rho(v_k^\text{acc}) }_2 \leq \tau \norm{r(v_{k-1})}_2
\end{equation}
using the previous residual norm $\norm{r(v_{k-1})}_2$ and $\tau \in (0, 2]$.
This means that every evaluation of~\eqref{eq:relaxed_safeguarding} only requires the evaluation of $F_\rho(\vacc)$.
We note that in the majority of iterations the safeguarding check~\eqref{eq:relaxed_safeguarding} passes, in which case $F_\rho(\vacc)$ can be used in the next evaluation of AA.
Consequently, the safeguarding step only leads to additional operator evaluations if the candidate point is rejected.

Another source of instability is the tendency of the iterates $r(v_k)$ to become nearly co-aligned, which tends to happen in regions where the underlying algorithm does not make much progress.
A consequence is bad conditioning of the matrix $\mc{R}_k$, which then results in unusable estimates of the coefficients $\eta_k$.
To avoid this scenario we monitor the norm of the coefficients.
If $\norm{\eta_k}_2 > \eta_{\max}$ we abort the acceleration step and perform an ordinary fixed-point iteration.

\subsection{Scheduling of operator changes}
Solver packages rarely use the same operator at each iteration but instead employ heuristics to tune parameters during the solve process.
For example, the ADMM step size parameter is often adapted based on the ratio of primal and dual residuals at the current iterate~\cite{Boyd_2011}.
However, AA relies on the history of past input vectors and operator outputs to approximate the inverse of its Jacobian and therefore requires a non-changing operator.
This history is invalid if the operator changes, and the method must then be restarted.
Another issue is that some solver packages rely on pure operator steps to monitor convergence behaviour.
For example~\cite{Banjac_2019} uses successive differences in FOM iterates and separating hyperplane conditions to detect infeasible problems.

We accommodate these two requirements by using a restarted memory approach for AA which adds information of past iterates up to a maximum memory length $m_{\max}$.
Then the memory is reset and begins building a new history starting with the last iterate.
Operator changes, i.e.\ computing a new $F_\rho$ based on a parameter update rule $\rho^{*} \leftarrow u(F_{\rho}, v_k, f_k, r_k, \rho)$, are scheduled after the AA method has been restarted or the safeguarding check~\eqref{eq:relaxed_safeguarding} fails, so that any changes in the operator are taken into account when new iterates are collected.
Similarly, the infeasibility detection is scheduled after AA restarts, because each restart is followed by at least two non-accelerated iterations.

\section{NUMERICAL RESULTS}
\label{sec:benchmarks}
To evaluate the impact of safeguarding and acceleration on a FOM, we implemented Algorithm~\ref{alg:safeguarded_acceleration} in v$0.8$ of the conic ADMM-based solver \pkg{COSMO}.
Different variants of AA are implemented in a standalone package\footnote{\href{https://github.com/oxfordcontrol/COSMOAccelerators.jl}{{https://github.com/oxfordcontrol/COSMOAccelerators.jl}}}.
Here we compare three different configurations of \pkg{COSMO}: without acceleration, with \emph{unsafe} acceleration, and the safeguarded acceleration detailed in Algorithm~\ref{alg:safeguarded_acceleration}.
The experiments were run using Julia v$1.5$ on computing nodes of the University of Oxford ARC-HTC cluster with 16 logical Intel Xeon E5-2560 cores and 64GB of DDR3 RAM.
We used the default parameters of \pkg{COSMO} and set the accuracy to $\epsilon = \num{e-6}$ for QPs and $\epsilon = \num{e-5}$ for SDPs, checking for convergence every $25$ iterations.
The acceleration method was configured with maximum memory length $m_{\text{max}} = 15$.
This value was chosen to obtain the best trade-off between convergence benefits, singularity issues, and computation overhead of higher memory lengths.
For the safeguarding parameter we chose $\tau = \num{2}$, and for the maximum norm of the Anderson parameters $\eta_{\max} = \num{e4}$, as these values work reasonably well on many different problem types.

Our benchmark tests used more than 500 problems from seven different QP and SDP problem sets to evaluate the impact of acceleration on the number of iterations needed to achieve higher accuracy solutions.
In each case we evaluated how much extra computation time the acceleration scheme required.
The tests compare three QP problem sets: the Maros and M\'{e}sz\'{a}ros problem set~\cite{Maros_1999}; model predictive control problems based on the MPC Benchmarking Collection~\cite{Ferreau_2020}; and Markowitz portfolio optimisation problems of the form described in~\cite{Stellato_2020}.
Furthermore, problems from four SDP problem sets are included: a relaxed sparse principal component analysis (SPCA) problem~\cite{d_aspremont_2004}; computation of the Lov\'{a}sz theta function~\cite{Lovasz_1979} for a set of undirected graphs from the SuiteSparse Matrix collection~\cite{Davis_2015}; randomly generated SDPs with chordal block arrowhead sparsity pattern~\cite{Garstka_2020a}; and a set of smaller non-decomposable problems used in Hans Mittelmann's SDP benchmarks~\cite{Mittelmann_}.

To compare the three solver methods, we calculate for each problem set the mean and median total solve time, the number of problems solved, and the mean time used to calculate the Anderson candidate points as a fraction of the solve time.
The results are shown in Table~\ref{tab:acceleration_results}.

% Approach from: https://tex.stackexchange.com/questions/89115/how-to-rotate-text-in-multirow-table
\begin{table*}[htb]
\begin{center}
   \sisetup{round-mode=places
           ,round-precision=1
           ,scientific-notation=false}
\caption{Results for vanilla, accelerated, and safeguarded \& accelerated ADMM for various QP and SDP problem sets.}\label{tab:acceleration_results}
\begin{threeparttable}
% \begin{footnotesize}

 % \resizebox{0.8\textwidth}{!}{%
   \begin{tabular}{ c | l c l c c c}
   \toprule
     & algorithm & solved & iter$^1$ & solve time$^2$  &  \% acc time$^3$ &  gmean$^4$\\
     \midrule
     \parbox[t]{2mm}{\multirow{3}{*}{\rotatebox[origin=c]{90}{Maros}}}  & vanilla  & 75 & \num{1046.0144927536232}  & \num[round-precision=2]{2.1005679386249487} (\num[round-precision=2]{0.017707109451293945}) & \num{} & \num[round-precision=2]{2.3338367818703003} \\
 & accelerated  & 98 & \num{676.8115942028985}  & \num[round-precision=2]{1.2651565593221914} (\num[round-precision=2]{0.024251937866210938}) & \num{24.694235050365993} & \num[round-precision=2]{1.0845888224127997} \\
 & safeguarded  & 100 & \num{505.7971014492754} (\num{22.89855072463768}) & \num[round-precision=2]{1.1576446173847585} (\num[round-precision=2]{0.02468085289001465}) & \num{24.986034701644368} & \num[round-precision=2]{1.0} \\
\midrule
\parbox[t]{2mm}{\multirow{3}{*}{\rotatebox[origin=c]{90}{MPC}}}  & vanilla  & 255 & \num{820.273631840796}  & \num[round-precision=3]{0.019628104878895318} (\num[round-precision=4]{0.0026268959045410156}) & \num{} & \num[round-precision=2]{2.744748554966579} \\
 & accelerated  & 230 & \num{296.8905472636816}  & \num[round-precision=3]{0.01247470414460595} (\num[round-precision=4]{0.0024809837341308594}) & \num{27.500246500378257} & \num[round-precision=2]{4.717986231250305} \\
 & safeguarded  & 285 & \num{272.0149253731343} (\num{115.28855721393035}) & \num[round-precision=3]{0.014768876839633011} (\num[round-precision=4]{0.0028579235076904297}) & \num{26.46294630694616} & \num[round-precision=2]{1.0} \\
\midrule
\parbox[t]{2mm}{\multirow{3}{*}{\rotatebox[origin=c]{90}{Portfolio}}}  & vanilla  & 10 & \num{3460.0}  & \num[round-precision=2]{171.50901823043824} (\num[round-precision=2]{134.16152358055115}) & \num{} & \num[round-precision=2]{3.041892313078642} \\
 & accelerated  & 10 & \num{1042.5}  & \num[round-precision=2]{69.85312557220459} (\num[round-precision=2]{48.96264696121216}) & \num{4.7673069307595215} & \num[round-precision=2]{1.31122936015046} \\
 & safeguarded  & 10 & \num{725.0} (\num{30.7}) & \num[round-precision=2]{49.6727395772934} (\num[round-precision=2]{32.292999505996704}) & \num{3.648461311759099} & \num[round-precision=2]{1.0} \\
\midrule
\parbox[t]{2mm}{\multirow{3}{*}{\rotatebox[origin=c]{90}{SPCA}}}  & vanilla  & 9 & \num{5586.111111111111}  & \num[round-precision=2]{132.81074661678738} (\num[round-precision=2]{95.69017195701599}) & \num{} & \num[round-precision=2]{6.724995205775979} \\
 & accelerated  & 10 & \num{1013.8888888888889}  & \num[round-precision=2]{46.917398585213554} (\num[round-precision=2]{15.597301006317139}) & \num{14.57337557048512} & \num[round-precision=2]{1.571301830376615} \\
 & safeguarded  & 10 & \num{652.7777777777778} (\num{32.333333333333336}) & \num[round-precision=2]{23.67362485991584} (\num[round-precision=2]{12.203541994094849}) & \num{13.513839932479025} & \num[round-precision=2]{1.0} \\
\midrule
\parbox[t]{2mm}{\multirow{3}{*}{\rotatebox[origin=c]{90}{Block}}}  & vanilla  & 20 & \num{3136.25}  & \num[round-precision=2]{99.30660681724548} (\num[round-precision=2]{42.10937297344208}) & \num{} & \num[round-precision=2]{5.402246143620185} \\
 & accelerated  & 20 & \num{308.75}  & \num[round-precision=2]{11.488829147815704} (\num[round-precision=2]{8.465690016746521}) & \num{2.9406634886926497} & \num[round-precision=2]{1.0} \\
 & safeguarded  & 20 & \num{396.25} (\num{66.9}) & \num[round-precision=2]{18.847397172451018} (\num[round-precision=2]{8.966078042984009}) & \num{2.2580617775537317} & \num[round-precision=2]{1.2654263861928052} \\
\midrule
\parbox[t]{2mm}{\multirow{3}{*}{\rotatebox[origin=c]{90}{Lov\'{a}sz}}}  & vanilla  & 12 & \num{2220.8333333333335}  & \num[round-precision=2]{17.012833495934803} (\num[round-precision=2]{7.21196448802948}) & \num{} & \num[round-precision=2]{5.72367439636774} \\
 & accelerated  & 15 & \num{1360.4166666666667}  & \num[round-precision=2]{5.787716845671336} (\num[round-precision=2]{3.0127739906311035}) & \num{8.748547678621183} & \num[round-precision=2]{1.091143824279882} \\
 & safeguarded  & 15 & \num{1283.3333333333333} (\num{15.0}) & \num[round-precision=2]{4.922293384869893} (\num[round-precision=2]{2.817574381828308}) & \num{8.544940221229872} & \num[round-precision=2]{1.0} \\
\midrule
\parbox[t]{2mm}{\multirow{3}{*}{\rotatebox[origin=c]{90}{Mittelm.}}}  & vanilla  & 22 & \num{1853.409090909091}  & \num[round-precision=2]{179.50180557641116} (\num[round-precision=2]{34.81064462661743}) & \num{} & \num[round-precision=2]{1.716933904618074} \\
 & accelerated  & 29 & \num{744.3181818181819}  & \num[round-precision=2]{79.56544262712652} (\num[round-precision=2]{31.836288332939148}) & \num{3.675362209283791} & \num[round-precision=2]{1.0} \\
 & safeguarded  & 31 & \num{792.0454545454545} (\num{5.909090909090909}) & \num[round-precision=2]{88.68485282767902} (\num[round-precision=2]{22.51495897769928}) & \num{3.4259432075556076} & \num[round-precision=2]{1.0090829344515913}
     \\ % need an extra linebreak here to avoid compilation errors
     \bottomrule
   \end{tabular}
   % }
  % \end{footnotesize}

 \begin{tablenotes}
 \scriptsize{
     \item[1] mean iteration and (extra safeguarding iterations);\\
     \item[2] mean and median (based on subset of problems where all solver configurations solved the problem);\\
    \item[3] geometric mean of fraction of total solve time spent in acceleration-related functions; \item[4] normalized shifted geometric mean of solve time, see~\eqref{eq:shifted_gmean} (based on all problems in the problem set);}
 \end{tablenotes}

 \end{threeparttable}
\end{center}
\end{table*}

As each solver configuration solved a different number of problems, we determined the average iteration counts and the average solve times using the subset of problems that was solved by every solver configuration.
This skews the results slightly in favour of solver configurations that solved fewer problems as problems not solved by each configuration tend to converge slower.
For the safeguarded method, Table~\ref{tab:acceleration_results} also shows in brackets the number of additional operator evaluations due to declined candidate points that failed the safeguarding check~\eqref{eq:relaxed_safeguarding}.
Another metric shown in Table~\ref{tab:acceleration_results}, and commonly used in solver benchmarks, is the normalized shifted geometric mean $\mu_{g, s}$ of the solve time, defined by
\begin{equation}\label{eq:shifted_gmean}
%  \mu_{g, s} \eqdef \sqrt[n]{\prod_{p}(t_{p, s} + \mathrm{sh}) - \mathrm{sh} }
  \mu_{g, s} \eqdef \Bigl[ \prod_{p}(t_{p, s} + \mathrm{sh}) - \mathrm{sh} \Bigr]^{1/n}
\end{equation}
with total solver time $t_{p, s}$ of solver $s$ and problem $p$, shifting factor $\mathrm{sh}$ and size of the problem set $n$.
We used a shifting factor $\mathrm{sh}=10$ and a maximum allowable time of \SI{5}{\min} for QPs and \SI{60}{\min} for SDPs.
Compared to the vanilla (i.e.\ unaccelerated) method, both accelerated methods provide a significant reduction in both the mean number of iterations and solve time.
For the Maros and MPC problem sets, the median solve time stays fairly similar to the vanilla method. This is due to the presence of a number of easy problems that are solved in only a few iterations by each method.

The results also show that the safeguarded acceleration leads to a greater number of problems solved.
While the impact of the safeguarded vs.\ non-safeguarded method appears small for most problem sets, we see a much more robust behavior for the MPC problems, with an additional 55 problems solved.
The different numbers of problems solved is taken into account by the shifted geometric mean of the solve time.
The impact of the safeguarded acceleration vs.\ the vanilla method ranges between $1.72$ for the Mittelmann SDPs and $6.72$ for the SPCA problems.
To evaluate the additional time spent inside acceleration related functions, we compute the time spent on acceleration as a fraction of the total solve time.
For a memory size of \num{15} this fraction varies from \SIrange{2}{15}{\percent} for the large problems, the Markowitz Portfolio and SDP problems.
For the smaller QP problem sets this fraction grows to \SI{25}{\percent}.
Consequently, for smaller problems a significant reduction in iterations is needed to reduce the overall solve time.
On average this seems to be the case for the Maros and MPC problem sets, but for some problems the acceleration might slow the solver down.

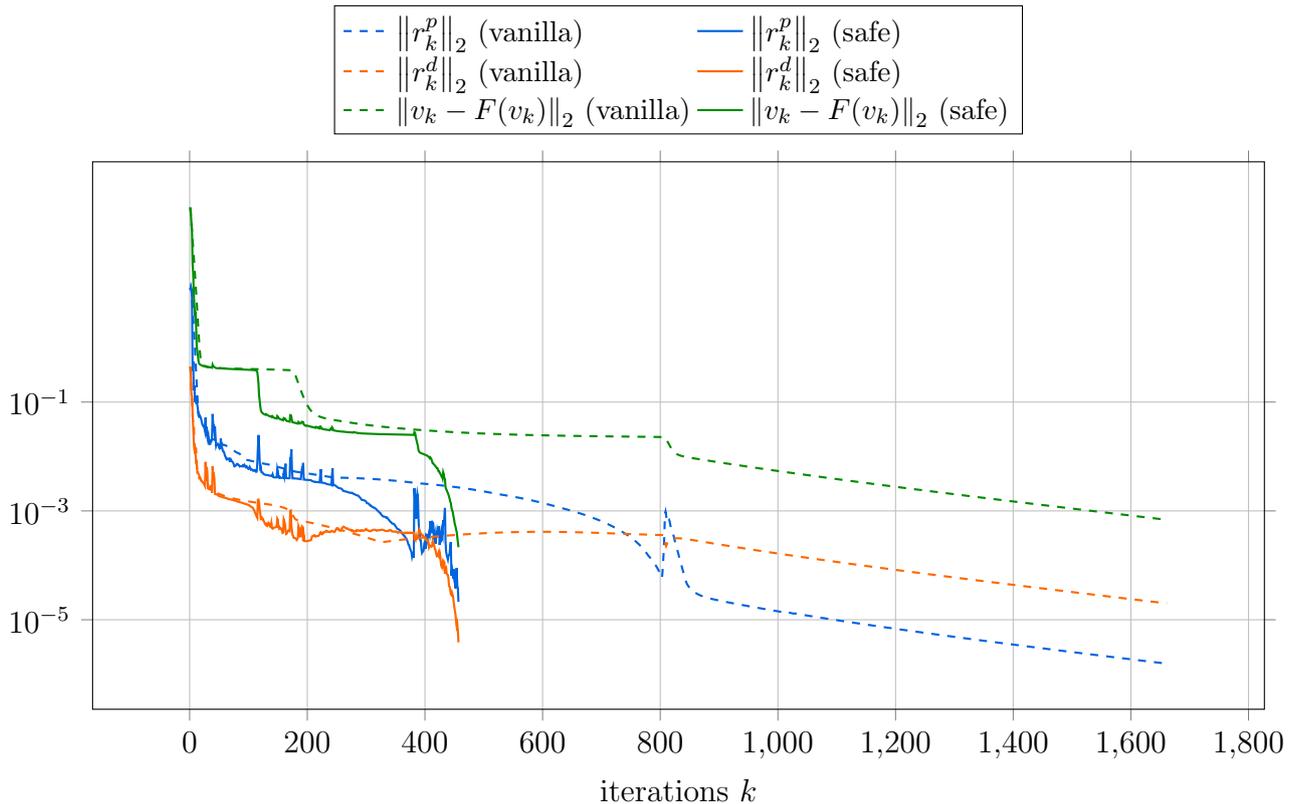
\begin{figure}[htb]
  \centering
  \vspace{5mm}
  \begin{tikzpicture}
  \begin{semilogyaxis}[
  width=\linewidth,
  height=.52\linewidth,
  xlabel={iterations $k$},
  ytick={0.1,0.001,0.00001},
  tick align=outside,
  grid=both,
  minor y tick num=10,
  % y grid style={lightgray!92.026143790849673!black},
  legend style={at={(0.5,1.05)},anchor=south},
  legend cell align={left},
  legend columns=2
  ]
    \addplot [thick, dashed, cosmocolor] table [x={iter}, y={vanillaPrimal}]{./residual_plot_ros_500_vanilla.csv};
      \addlegendentry{\small{$\norm{r_k^p}_2$ (vanilla)}}
    \addplot [thick, cosmocolor] table [x={iter}, y={safePrimal}]{./residual_plot_ros_500_safeguarded.csv};
      \addlegendentry{\small{$\norm{r_k^p}_2$ (safe)}}
    \addplot [thick, dashed, cosmocdcolor] table [x={iter}, y={vanillaDual}]{./residual_plot_ros_500_vanilla.csv};
     \addlegendentry{\small{$\norm{r_k^d}_2$ (vanilla)}}
    \addplot [thick, cosmocdcolor] table [x={iter}, y={safeDual}]{./residual_plot_ros_500_safeguarded.csv};
     \addlegendentry{\small{$\norm{r_k^d}_2$ (safe)}}
    \addplot [thick, dashed, scscolor] table [x={iter}, y={vanillaFPNorm}]{./residual_plot_ros_500_vanilla.csv};
      \addlegendentry{\small{$\norm{v_k - F(v_k)}_2$ (vanilla)}}
    \addplot [thick, scscolor] table [x={iter}, y={safeFPNorm}]{./residual_plot_ros_500_safeguarded.csv};
      \addlegendentry{\small{$\norm{v_k - F(v_k)}_2$ (safe)}}
  \end{semilogyaxis}

\end{tikzpicture}
   \caption{Norms of primal residual $\norm{r_p^k}_2$, dual residual $\norm{r_d^k}_2$, and fixed point residual $\norm{v_k - F_\rho(v_k)}_2$ of the vanilla and the safeguarded accelerated method (Mittelmann: \texttt{ros\_500}).}
  \label{fig:residuals_ros500}
\end{figure}

Figure~\ref{fig:residuals_ros500} shows the convergence of the ADMM and operator residuals for the vanilla and the safeguarded accelerated method for the SDP problem \texttt{ros\_500}.
This behavior is typical when AA works successfully.
Initially the residuals of both methods decrease in a similar fashion until an accuracy of $10^{-2}$ or $10^{-4}$ is reached.
Then the residuals of the accelerated method (solid lines) drop sharply relative to the vanilla method's
%and avoid the vanilla method's period of
slow convergence from iteration 400 onwards.
This is likely due to the algorithm reaching the region where AA approximates the Jacobian well and
%the method
achieves superlinear convergence.
It also indicates that the impact of AA will be fairly low when used to accelerate the convergence only to a low accuracy of e.g. $10^{-3}$.

\section{CONCLUSIONS}
\label{sec:conclusion}
This paper uses a combination of scheduled memory restarts, a safeguarding rule based on the residual operator norm, and least-squares condition checking, to safeguard Anderson acceleration.
We show that the approach works well for a FOM-based solver that allows adaptation of its operator form to improve convergence behaviour and which relies on successive unaltered iterates for infeasibility detection.
We provide an efficient  AA implementation using an updated QR decomposition in the latest version of \pkg{COSMO}.
The effectiveness of our approach in reducing both the mean number of iterations and the solve time while increasing the number of solved problems is shown for a large number of QPs and SDPs taken from different application domains.
Instead of using $m_k + 1$ past iterates in the acceleration scheme, it seems promising to investigate whether performance improvements are achievable using the same number of iterates but spread out over a longer history of past iterations.

\bibliographystyle{IEEEtran}
\bibliography{\rootPath/acceleration_references.bib}

\end{document}